\documentclass[lettersize,journal]{IEEEtran}

\usepackage{amsmath,amsfonts}
\usepackage{algorithm}
\usepackage{algpseudocode}
\usepackage{array}
\usepackage[caption=false,font=normalsize,labelfont=sf,textfont=sf]{subfig}
\usepackage{textcomp}
\usepackage{stfloats}
\usepackage{url}
\usepackage{verbatim}
\usepackage{graphicx}
\usepackage{cite}
\usepackage{mathtools}
\usepackage{amsthm}
\usepackage{orcidlink}
\usepackage{placeins}

\theoremstyle{definition}

\newcommand{\TP}{{d}}

\hypersetup{
    colorlinks = true,
    allcolors = {blue}
}

\begin{document}
\bstctlcite{bstctl:etal, bstctl:nodash, bstctl:simpurl}
\title{Towards Efficient Aggregation of Storage Flexibilities in Power Grids}

\author{Emrah Öztürk \orcidlink{0000-0001-5063-531X}, 
        Kevin Kaspar \orcidlink{0009-0007-6176-2737}, 
        Timm Faulwasser \orcidlink{0000-0002-6892-7406}, 
        Karl Worthmann \orcidlink{0000-0002-1450-2373},
        Peter Kepplinger, \orcidlink{0000-0003-2440-7270} \\ and
        Klaus Rheinberger \orcidlink{0000-0002-8277-4944}

\thanks{EÖ, KK, PK and KR: 
This work was supported by the project ”EBusCharge” [FFG, No. 899915]. This project is supported with funds from the Climate and Energy Fund and implemented in line with the "Zero Emission Mobility" programme. This work was also supported by the project "FreeE-Bus" (interreg, ABH006).}
\thanks{EÖ, KK, PK and KR are with the Research Center Energy, Vorarlberg University of Applied Sciences, illwerke vkw Endowed Professorship for Energy Efficiency, 6850 Dornbirn, Austria (e-mail: emrah.oeztuerk@fhv.at; kevin.kaspar@fhv.at; peter.kepplinger@fhv.at; klaus.rheinberger@fhv.at).}
\thanks{TF is with the Institute of Energy Systems, Energy
Efficiency and Energy Economics, TU Dortmund University, 44227 Dortmund, Germany (e-mail: timm.faulwasser@ieee.org).}
\thanks{KW is with the Institute of Mathematics,
Technische Universität Ilmenau, 98693 Ilmenau, Germany (e-mail:
karl.worthmann@tu-ilmenau.de).}}

\maketitle

\begin{abstract}
The increasing penetration of volatile renewables combined with increasing demands poses a challenge to modern power grids. Furthermore, distributed energy resources and flexible devices (electric vehicles, PV generation, ...)  are becoming more widespread, making their aggregate usage for ancillary services interesting. However, accurately quantifying the aggregate flexibility of numerous flexible devices is known to be limited by the curse of dimensionality, i.e., it does not scale well computationally. This has led to the development of various approximation algorithms. In this study, we improve upon our previously proposed vertex-based inner approximation, extending it to more general storage devices. We demonstrate the efficacy and accuracy of the proposed method in a case study comparing our approach with an exact centralized control framework, where the flexibility of numerous electric vehicles is combined to reduce the peak load in a residential area.
\end{abstract}


\section{Introduction}\label{sec:introduction}
As the world moves towards sustainable energy sources, the use of renewables such as solar and wind power has become widespread. 
Simultaneously, there has been a notable increase in electricity consumption. The inherent variability of renewable energy generation combined with increasing electricity demand poses a significant risk to the reliable and secure operation of modern power grids. 
At the same time, flexible devices, such as stationary batteries, electric vehicles (EVs), thermostatically controlled loads, etc., are on the rise and offer flexibility potential. 
By effectively managing these devices, it is possible to shift electricity demand or to store excess generation for later use during periods of high demand. 
In this sense, numerous flexible devices may act as a virtual power plant. In response to the need for a framework that allows large-scale control, aggregators are proposed as a viable concept \cite{gkatzikis_role_2013}. 
Aggregators manage contracted user devices by quantifying their overall flexibility, by participating in electricity markets, and by distributing power profiles to individual devices. 
This enables aggregators to offer various services such as demand response, frequency regulation, or energy arbitrage at scale. 
In return for participating in this framework device owners receive revenues. However, the pivotal bottleneck in aggregation-based control and coordination is the numerical scalability of the aggregation problem, i.e., the exact quantification of the aggregate flexibility. 
As the time horizon grows and as the number of devices increases, this problem becomes computationally intractable \cite{tiwary_hardness_2008}. 
For this reason, several approximate approaches have been proposed in the literature. Existing approximation methods can be roughly divided into inner and outer approximations. 
Inner approximations \cite{ozturk_alleviating_2024, zhao_geometric_2017, hao_aggregate_2015, he_hao_generalized_2013, nazir_inner_2018, muller_aggregation_2015, muller_aggregation_2019, zhen_computing_2018, appino_towards_2021, al_taha_multi-battery_2023, zhang_fast_2024,jian_analytical_2024} are conservative, i.e., the feasible region is potentially reduced, excluding feasible power profiles. 
Outer approximations \cite{zhao_geometric_2017, hao_aggregate_2015, he_hao_generalized_2013, zhang_fast_2024, barot_outer_2016, barot_concise_2017} on the other hand, overestimate the feasible region and, thus, may contain infeasible power profiles. 
Auxiliary service providers must ensure the fulfilment of the power profiles to avoid contractual penalties, which favors inner approximations. 
However, state-of-the-art inner approximations suffer from objective-dependent performance, high computational complexity, and might exclude nominal power profiles such as the vector of zeros, cf. \cite{ozturk_aggregation_2022}. 
In our previous study \cite{ozturk_alleviating_2024}, we developed a vertex-based inner approximation that mitigates these problems. 
Our approach outperformed ten state-of-the-art inner approximations in the open-source aggregation benchmark \cite{ozturk_aggregation_2022}. 
However, there are still limitations with respect to the applicability of this method. 
Hence, in this paper we extend the method from \cite{ozturk_alleviating_2024} to cover a wider range of storage models.

The contributions can be summarized as follows: Building upon our previous research \cite{ozturk_alleviating_2024}, we propose a new algorithm that can handle various storage models with limited availability and varying energy constraints. 
We demonstrate the applicability of the approach in a case study and show its effectiveness in realistic scenarios. 

The remainder of this paper is organized as follows: In Section \ref{sec:2}, the general storage model is introduced, and its applications are discussed. 
Section \ref{sec:3} presents the novel algorithm for quantifying the aggregated flexibility of a fleet of storage devices. 
Case studies are presented in Section \ref{sec:4}. Lastly, the paper ends with conclusions in Section \ref{sec:5}.

\section{Preliminaries}\label{sec:2}
In this section, we present a general model capable of representing a variety of energy storage devices. 
The applicability of this model is discussed for specific energy-related applications. 
We consider a discrete-time setting with $d$ time intervals each of duration $\Delta t$. 
Let $x_t$ be the (dis-)charging power during period $t$ and let $S_t$ the energy in the storage after period $t$. 
A general storage device is subject to the following constraints and system dynamics:
\begin{subequations}\label{eq:1}
\begin{align}
    &\underline{x}_t \leq x_t \leq \overline{x}_t \quad \forall t = 1, \ldots, d\\
    &S_t = \alpha S_{t-1} + x_t\Delta t \quad \forall t = 1, \ldots, d\\
    &\underline{S}_t \leq S_t \leq \overline{S}_t \quad \forall t = 0, \ldots, d\\
    &S_0 = S_\text{init}
\end{align}
\end{subequations}
where $\underline{x}_t,\  \overline{x}_t, \underline{S}_t,\ \overline{S}_t$ denote the lower and upper limits for the variables $x_t$ and $S_t$, respectively. 
$S_\text{init}$ denotes the initial energy and $\alpha \in (0,\ 1]$ the self-discharge factor. Eqs. \eqref{eq:1} define the polytope
\begin{align}
    &\mathcal{B}(\underline{x},\overline{x},\underline{S}, \overline{S},\alpha, S_\text{init}) \coloneq\\
    &\{x \in \mathbb{R}^\TP : A(\alpha)x \leq b(\underline{x},\overline{x},\underline{S}, \overline{S}, \alpha, S_\text{init})\},\nonumber
\end{align}    
where $A \in \mathbb{R}^{4\TP\times \TP}$ and $b \in \mathbb{R}^{4\TP}$ are defined by
\begin{subequations}\label{eq:Ab}
	\begin{align}
		&A(\alpha) \coloneqq \left(
		-I , I , \Gamma^\top , -\Gamma^\top
		\right)^\top \\
		&b(\underline{x},\overline{x},\underline{S}, \overline{S}, \alpha, S_\text{init}) \coloneqq\\ 
		&\left(-\underline{x}^\top, \overline{x}^\top, \frac{1}{\Delta t}(\overline{S} - S_\text{init}a_d)^\top, \frac{1}{\Delta t}(-\underline{S}+S_\text{init}a_d)^\top\right)^\top. \nonumber
	\end{align}
\end{subequations}
Moreover, we have $a_\TP \coloneqq (\alpha, \alpha^2, \ldots, \alpha^\TP)^\top$ and $I \in \mathbb{R}^{\TP \times \TP}$ is the identity matrix, and $\Gamma \in \mathbb{R}^{\TP \times \TP}$ is a Toeplitz matrix with first column and row defined by $(1, \alpha, \ldots, \alpha^{\TP-1})^\top$ and $(1, 0, \ldots, 0)$, respectively. 
The storage model \eqref{eq:1} can be used to represent various real-world appliances, like battery energy storage systems, EVs, and thermostatically controlled loads. 
For instance, an EV can be modelled by using the following bounds in \eqref{eq:Ab}
\begin{subequations}
\begin{align}
    & \overline{x}_t = x_\text{max}\lambda_t, \ \underline{x}_t = x_\text{min}\lambda_t, \quad \forall t = 1, \ldots, d,\\
    & \overline{S}_t = S_\text{max} + \sum_{\tau = 1}^t\alpha^{t-\tau}D_\tau\Delta t, \quad  \forall t = 1, \ldots, d,\\
    & \underline{S}_t = S_\text{min} + \sum_{\tau = 1}^t\alpha^{t-\tau}D_\tau\Delta t, \quad  \forall t = 1, \ldots, d-1,\\
    & \underline{S}_d = S_f + \sum_{\tau = 1}^d\alpha^{t-\tau}D_\tau\Delta t,
\end{align}
\end{subequations}
where $x_\text{max},x_\text{min}, S_\text{max}, S_\text{min} \in \mathbb{R} $ are the maximum and minimum power and energy bounds, respectively, $\lambda\in \{0,1\}^d$ is the availability vector, $D \in \mathbb{R}^d$ is the trip consumption vector, and $S_f \in \mathbb{R}$ the minimum final energy. 
A stationary battery can be represented by the special case of an EV, where trip consumption is the vector of zeros and the availability the vector of ones.

\section{Aggregation of Storage Devices}\label{sec:3}
We present a novel algorithm for quantifying the aggregate flexibility of a fleet of storage devices of the type $\mathcal{B}\left(\underline{x},\overline{x},\underline{S},\overline{S},\alpha,S_\text{init}\right)$. 
To this end, extreme actions, ideally vertices or elements located on the facets of the polytopes, are computed within the storage device $\mathcal{B}\left(\underline{x},\overline{x},\underline{S},\overline{S},\alpha, S_\text{init} \right)$. 
Corresponding extreme actions are then summed to give elements in the aggregate flexibility $\{x\in\mathbb{R}^d:x=\sum_{i=1}^{n}x_i, x_i \in \mathcal{B}\left({\underline{x}}_i,{\overline{x}}_i,{\underline{S}}_i,{\overline{S}}_i,\alpha_i, S_{\text{init},i} \right)\}$. 
The area spanned by the convex hull of the summed extreme actions finally results in an inner approximation. 

\begin{figure}[tb]
    \centering
    \includegraphics[width=0.5\textwidth]{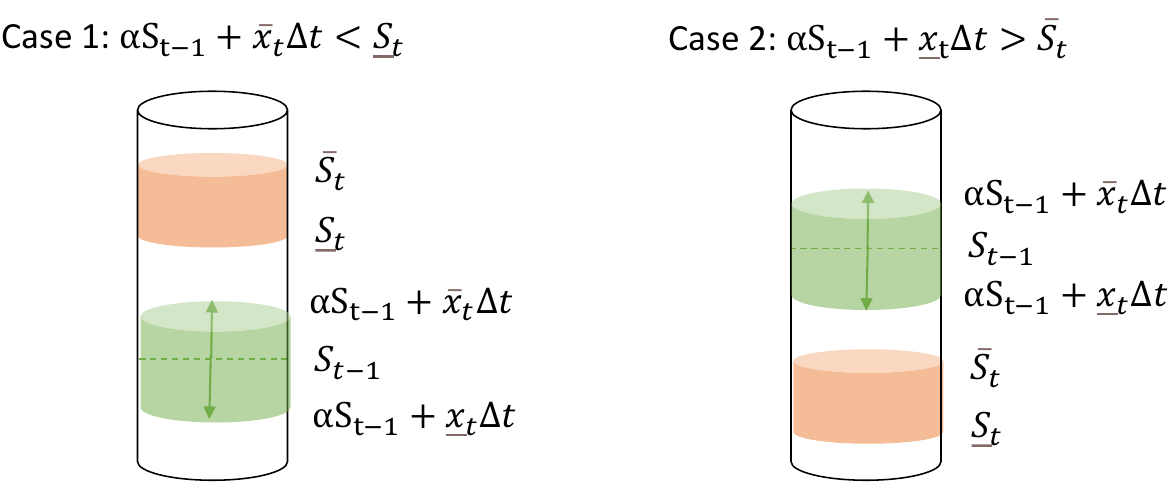}
    \caption{Example cases, in which the energy constraints in Eqs. \eqref{eq:1} are violated}
    \label{fig:violations}
\end{figure}

One approach to generate extreme actions is to move as far as possible in specific directions. 
To this end, we use $-1$ to indicate the negative direction and $1$ to indicate the positive direction. 
A direction in $d$-dimensional space can then be represented by a vector $j\in\{-1,1\}^d$. 
\begin{algorithm}[tb]
	\caption{(extremeActions)}\label{algo:extreme action}
	\begin{algorithmic}[1]
        \Require $\underline{x}, \overline{x}, \underline{S}, \overline{S}, S_\text{init}, \alpha, \mathcal{J} \subseteq \{-1, 1\}^\TP$
        \Ensure $V$
        \State $V \gets \mathbf{0}_{\TP \times |\mathcal{J}|}$ \Comment{$\TP \times |\mathcal{J}|$ matrix of zeros}
        \For{$j \in \mathcal{J}$}
            \State $k \gets 1$
            \State $y^j \gets \mathbf{0}_\TP$
            \For{$t = 1$ \textbf{to} $\TP$}
                \If{$j_t = 1$}
                    \State $p \gets \frac{\overline{S}_t - (\alpha^t S_\text{init} + \sum_{\tau = 1}^{t-1}\alpha^{t-\tau}y^j_\tau \Delta t)}{\Delta t}$
                    \State $y^j_t \gets \max(\min(\overline{x}_t, p), \underline{x}_t)$
                    \State $y^j \gets$ correctiveIncrease($y^j, \underline{x}, \overline{x}, \underline{S}, \overline{S}, S_\text{init}, \alpha, t$)
                    \State $y^j \gets$ correctiveDecrease($y^j, \underline{x}, \overline{x}, \underline{S}, \overline{S}, S_\text{init}, \alpha, t$)
                \ElsIf{$j_t = -1$}
                    \State $p \gets \frac{\underline{S}_t-(\alpha^t S_\text{init} + \sum_{\tau = 1}^{t-1}\alpha^{t-\tau}y^j_\tau \Delta t)}{\Delta t}$
                    \State $y^j_t \gets \min(\max(\underline{x}_t, p), \overline{x}_t)$
                    \State $y^j \gets$ correctiveIncrease($y^j, \underline{x}, \overline{x}, \underline{S}, \overline{S}, S_\text{init}, \alpha, t$)
                    \State $y^j \gets$ correctiveDecrease($y^j, \underline{x}, \overline{x}, \underline{S}, \overline{S}, S_\text{init}, \alpha, t$)
                \EndIf
            \EndFor
            \State $V[:, k] \gets y^j_i$
            \State $k \gets k + 1$
        \EndFor
    \end{algorithmic}
\end{algorithm}
The Algorithm \ref{algo:extreme action} computes elements $y_t^j$, given storage parameters and a set of vectors $\mathcal{J}\subseteq\{-1,1\}^d$. 
This is done by moving as far as possible in directions $j_t$ for each time period t, i.e., charging or discharging the storage to the limits. 
Lines 8 and 13 enforce the power constraints, i.e., $\underline{x}_t \leq y_t^j \leq \overline{x}_t$. 
However, as the energy limits are allowed to vary over time, future energy constraints could be violated. 
Two cases are distinguished as illustrated in Fig. \ref{fig:violations}. 
\begin{algorithm}[ht]
    \caption{(correctiveIncrease)}\label{algo:positive correction batt}
	\begin{algorithmic}[1]
        \Require $y^j, \underline{x}, \overline{x}, \underline{S}, \overline{S}, S_\text{init}, \alpha, t$
        \Ensure $y^j$
        \If{$\underline{S}_t > \alpha^t S_\text{init} + \sum_{\tau = 1}^{t} \alpha^{t-\tau}y_{\tau}^{j}\Delta t$}
		\State init $i$ \Comment{first index with $\overline{x}_i > 0$ in reversed $\overline{x}_{[t]}$}
        \State $p \gets \frac{\underline{S}_t - (\alpha^t S_\text{init} + \sum_{\tau = 1}^{t-1} \alpha^{t-\tau}y^j_{ \tau}\Delta t)}{\Delta t}$
        \State $y_{t-i+1}^{j} \gets \min(\max(\underline{x}_{t-i+1}, p), \overline{x}_{t+i+1})$
		\State $k \gets t - i$
		\While{$\underline{S}_t \neq \alpha^t S_\text{init} + \sum_{\tau = 1}^{t}\alpha^{t -\tau} y^j_{ \tau}\Delta t$}
        \For{$l=k$ to $t-i$}
            \State $p \gets \frac{\overline{S}_l - (\alpha^l S_\text{init} + \sum_{\tau = 1}^{l-1} \alpha^{l-\tau}y^j_{ \tau}\Delta t)}{\Delta t}$
            \State $y_{l}^j \gets \max(\min(\overline{x}_l, p), \underline{x}_l)$
        \EndFor
        \State $p \gets \frac{\underline{S}_t - (\alpha^t S_\text{init} + \sum_{\tau = 1}^{t-1} \alpha^{t-\tau}y^j_{ \tau}\Delta t)}{\Delta t}$
        \State $y_{t-i+1}^{j} \gets \min(\max(\underline{x}_{t-i+1}, p), \overline{x}_{t+i+1})$
		\State $k \gets k - 1$
		\EndWhile
		\EndIf
	\end{algorithmic}
\end{algorithm}

The case where the energy in the storage is less than the lower bound, left in Fig. \ref{fig:violations}, is handled via the corrective algorithm to increase the energy, cf. Algorithm \ref{algo:positive correction batt}. In Line 2, $i$ is initialized with the first index of reversed $\overline{x}_{\left[t\right]}$, where $\overline{x}_i>0$ applies. 
Here $\overline{x}_{\left[t\right]}\in\mathbb{R}^t$ refers to the vector consisting of the first $t$ elements of $\overline{x}$. 
This step is necessary to identify the last index prior that allows for correction, i.e., positive power. 
Correction to increase the energy must be carried out with strictly positive values and is not possible with indices $l>t-i+1$, as they are either zero or negative due to $\overline{x}_l \leq 0$. 
The power profile is then changed backwards without violating the constraints, starting with the index $t-i+1$ to reach $\underline{S}_t$. 
Lines 4, 9, and 12 enforce the power constraints, i.e., it holds that $\underline{x}_{t-i+1} \leq y_{t-i+1}^j \leq \overline{x}_{t-i+1}$. 
In Line 9, the $y_k^j$ are modified to increase the energy in the storage. 
Lines 4, and 12 guarantee that $\underline{S}_t$ is reached exactly. 
The algorithm converges for non-empty polytopes. 
In the limiting case, $y_k^j$ is set to the maximum possible value in each time period, and if $\underline{S}_t$ cannot be reached with maximum charging in each time period, then the set of feasible power profiles is empty.

The second case, right in Fig. \ref{fig:violations}, where the energy in the storage is greater than $\overline{S}_t$, is covered similarly in the corrective algorithm to decrease the energy, cf. Algorithm \ref{algo:negative correction batt}. The corrections, if needed, are applied in Lines 9, 10, 14, and 15 of Algorithm~\ref{algo:extreme action}.

 \begin{algorithm}[ht]
	\caption{(correctiveDecrease)}\label{algo:negative correction batt}
	\begin{algorithmic}[1]
        \Require $y^j, \underline{x}, \overline{x}, \underline{S}, \overline{S}, S_\text{init}, \alpha, t$
        \Ensure $y^j$
        \If{$\overline{S}_t < \alpha^t S_\text{init} + \sum_{\tau = 1}^{t} \alpha^{t-\tau}y_{\tau}^{j}\Delta t$}
		\State init $i$ \Comment{first index with $\underline{x}_i < 0$ in reversed $\underline{x}_{[t]}$}
        \State $p \gets \frac{\overline{S}_t - (\alpha^t S_\text{init} + \sum_{\tau = 1}^{t-1} \alpha^{t-\tau}y^j_{ \tau}\Delta t)}{\Delta t}$
        \State $y_{t-i+1}^{j} \gets \max(\min(\overline{x}_{t-i+1}, p), \underline{x}_{t+i+1})$
		\State $k \gets t - i$
		\While{$\underline{S}_t \neq \alpha^t S_\text{init} + \sum_{\tau = 1}^{t}\alpha^{t -\tau} y^j_{ \tau}\Delta t$}
        \For{$l = k$ to $t-i$}
            \State $p \gets \frac{\underline{S}_l - (\alpha^l S_\text{init} + \sum_{\tau = 1}^{l-1} \alpha^{l-\tau}y^j_{ \tau}\Delta t)}{\Delta t}$
            \State $y_k^j \gets \min(\max(\underline{x}_l, p), \overline{x}_l)$
        \EndFor
        \State $p \gets \frac{\overline{S}_t - (\alpha^t S_\text{init} + \sum_{\tau = 1}^{t-1} \alpha^{t-\tau}y^j_{ \tau}\Delta t)}{\Delta t}$
        \State $y_{t-i+1}^{j} \gets \max(\min(\overline{x}_{t-i+1}, p), \underline{x}_{t+i+1})$
		\State $k \gets k - 1$
		\EndWhile
		\EndIf
	\end{algorithmic}
\end{algorithm}

Given $n$ devices, we now assume that sets of extreme actions $\{y_i^j:j\in\mathcal{J}\}$ for $i=1,\ \ldots,\ n$ were computed using Algorithm \ref{algo:extreme action}. 
The sum of extreme actions with identical vector $j$ over all devices is denoted as:
\begin{equation}
    v^j \coloneq \sum_{i=1}^n y^j_i.
\end{equation}
By construction, the vector $v^j$ is an element of the aggregate flexibility. 
Hence, the convex hull of the set of summed vectors, i.e., $\text{Conv}\left(\{v^j:j\in\mathcal{J}\}\right)$ results in an inner approximation of the aggregate flexibility. 
The summed vectors can be constructed by adding the matrices given by Algorithm \ref{algo:extreme action}, i.e., $\sum_{i=1}^{n}V_i$ for devices $i\ =\ 1,\ \ldots,\ n$. 

\begin{figure}[bt]
    \centering
    \includegraphics[width=0.5\textwidth]{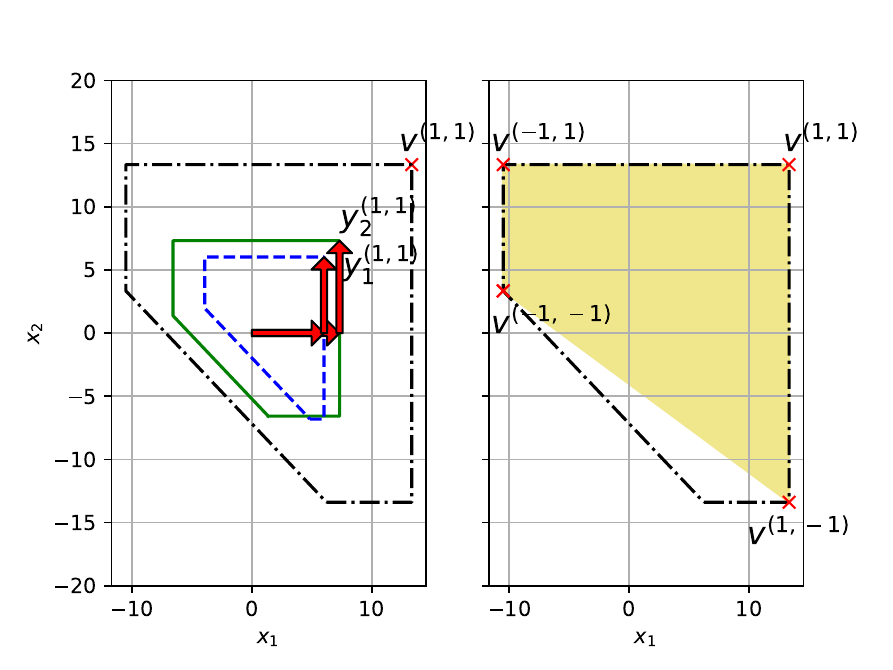}
    \caption{Left: extreme actions $y_1^{\left(1,1\right)},y_2^{\left(1,1\right)}$ in their polytopes, shown in dashed blue and solid green together with the summed extreme action $v^{\left(1,1\right)}$ in the Minkowski sum in dashed-dotted black. Right: all summed extreme actions $v^j, j \in \{-1,1\}^2$ with their convex hull.}
    \label{fig:operation}
\end{figure}
Since extreme actions are computed for every vector $j\in\mathcal{J}$ and for each device, the amount of extreme actions computations is given by $\left|\mathcal{J}\right|n$. 
For high-dimensional spaces, e.g., $d > 8$, we propose to select a set of vectors $\mathcal{J}$ by selecting $g$ distinct uniformly distributed elements in $\{-1,1\}^d$, where $g$ must be a function of the time periods, e.g., a hypercube has $2^d$ vertices in $d$-dimensional space. 
Fig. \ref{fig:operation} illustrates the operation of the proposed algorithm, where for the sake of simplicity, polytopes were used for which no correction is required. 
The extreme actions $y_1^j,y_2^j$ in their polytopes, together with their sum $v^j$ in the aggregated feasible region are shown on the left for the vector $(1, 1)$. 
Whereas on the right, all possible $v^j,j\in\{-1,1\}^2$ are shown together with their convex hull. 

Note that the convex hull is depicted for illustrative purposes only and does not need to be computed. 
The summed extreme actions are sufficient to describe the approximate feasible region and can be passed for optimization purposes. 
It should also be noted that for battery energy storage systems, the proposed approach simplifies to the method proposed in \cite{ozturk_alleviating_2024}. 
Therin, we have shown that the summed extreme actions are vertices of the aggregate feasible region.

\section{Case Study}\label{sec:4}
\begin{figure}[!b]
    \centering
    \includegraphics[width=0.4\textwidth]{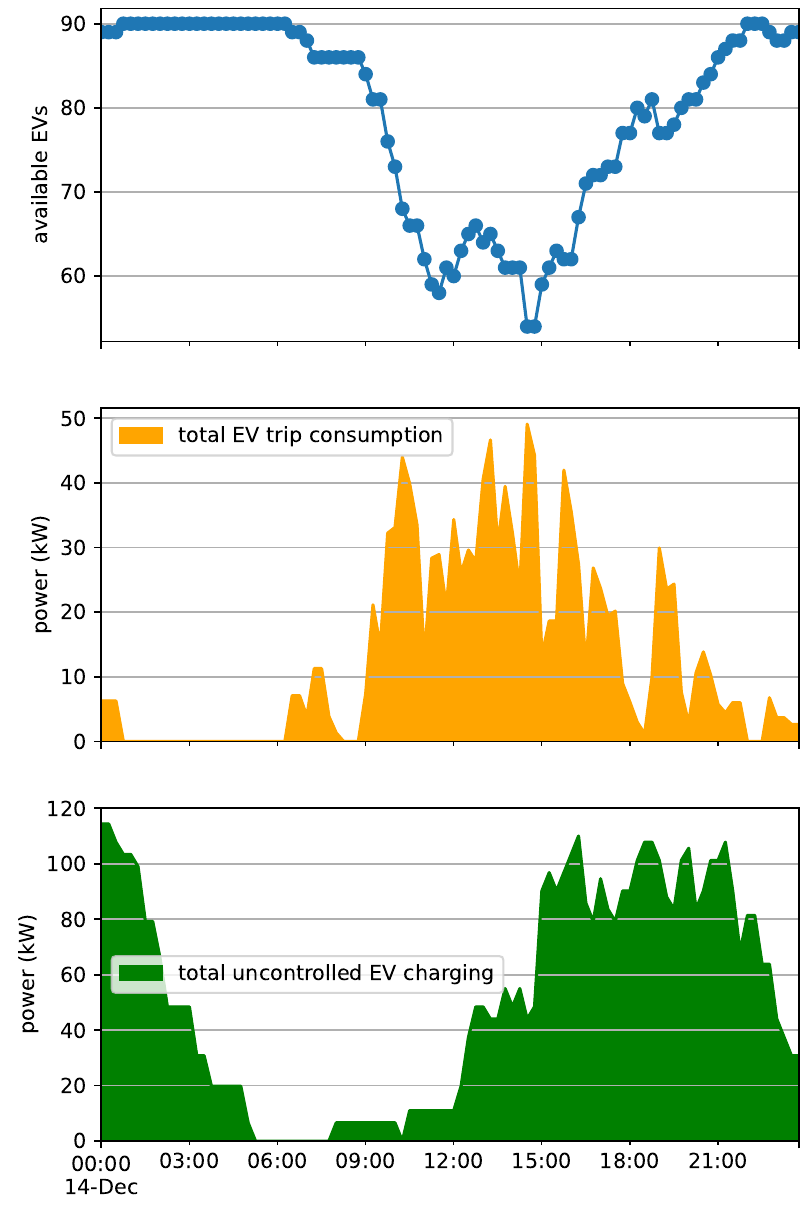}
    \caption{Available EVs (top), total EV trip consumption (centre), and total uncontrolled EV charging (bottom).}
    \label{fig:EV_data}
\end{figure}

\begin{figure}[tb]
    \centering
    \includegraphics[width=0.5\textwidth]{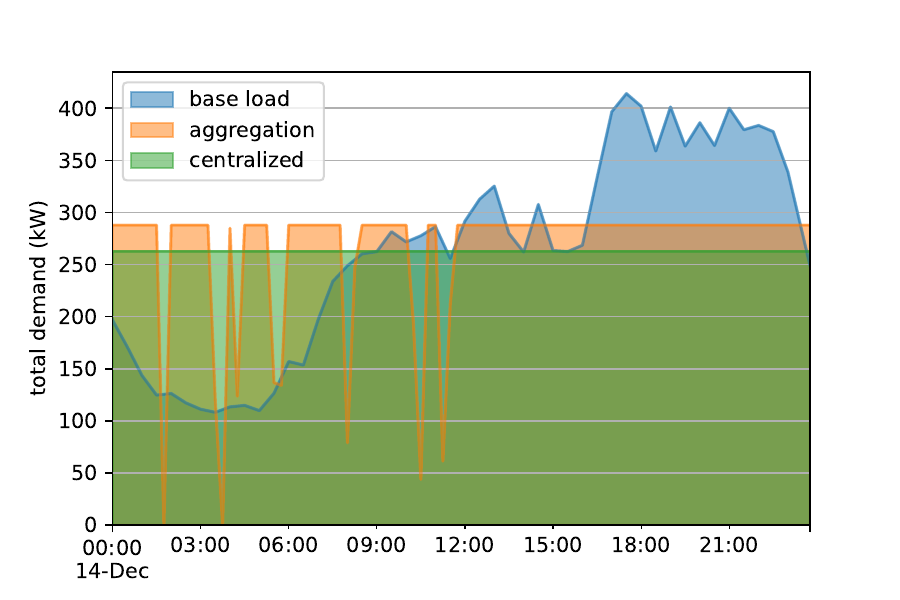}
    \caption{Residential loads: the base load is shown in blue, the aggregate controlled load in orange, and the centralized controlled load in green.}
    \label{fig:results}
\end{figure}

The proposed aggregation method is tested in a residential area with 300 households. 
An aggregator manages the charging and discharging of contracted customers' EVs. 
We use a one-day time horizon with a discrete time setting of 15-minute intervals. 
The aim is to utilise the aggregated flexibility to provide ancillary services to a grid operator, who reduces the overall peak load of the residential area. 
The method proposed is compared with a centralised control scheme to demonstrate the accuracy and computational efficiency. 
In contrast to the approach proposed, the centralised controller uses all constraints at once, without aggregation.
Gurobi \cite{gurobi_optimization_llc_gurobi_2024} is used to solve the resulting optimisation problem. 
We assume that 30 \% of the households use EVs with vehicle-to-grid technology. 
The data for the household demand and EVs, i.e., availability, charging, and trip consumption, cf. Fig. \ref{fig:EV_data}, are taken from \cite{rheinberger_dsm-data_2021}. 
We use the parameters of the Nissan Leaf, i.e., $x_{\mathrm{max}}=6.6$ kW, $x_{\mathrm{min}}=-6.6$ kW, and $S_{\mathrm{max}}=39$ kWh. 
At start, the EV batteries are charged to 50 \%. 
At the end of the operating window, the EVs must retain a minimum energy level, which is determined by subtracting the trip consumption from the initial energy and adding the energy charged given by the data. 
The polytopes are in general distinct due to the specific data associated with the EVs. 
We use a quadratic time dependence, i.e., $g = d^2 = 96^2$ distinct, uniformly selected vectors in $\{-1,1\}^d$ to calculate the extreme actions. 
Fig. \ref{fig:results} shows the total uncontrolled residential demand compared to the load curve achieved by the centralised, as well as the aggregation-based approach. 
The peak demand is effectively reduced to $283.08$ kW by the proposed method and to $262.68$ kW by the centralised approach. 
The calculation times are $35$ seconds and $5$ seconds for the proposed approach and the centralised approach, respectively. 
Note that although the centralised approach performs better, it is only used for comparison purposes and is not practical due to the lack of data security, and the increased communication effort for large-scale problems. 
Our approach, on the other hand, masks the data and reduces the communication effort significantly, i.e., only the aggregator communicates with the purchaser of the ancillary services. 
In addition, the proposed method can be implemented on a simple microcontroller, while the centralised approach requires special, possibly commercial, software.

\section{Conclusions}\label{sec:5}
In this study, we presented a novel inner approximation method that can be used to aggregate various storage models, like EVs, thermostatically controlled loads, and battery energy storage systems. 
The effectiveness and accuracy of the method are illustrated in a practical example where the aggregated flexibility of numerous EVs is utilized to reduce the peak load in a residential area. 
The proposed approach reduced the peak load effectively by 32\%, while the centralized controller achieved a peak power reduction of 37\%. 
The calculation times were 35 seconds for the pro-posed approach and 5 seconds for the centralized controller. 
Unlike the centralized approach, the aggregation method proposed masks the data, i.e., protects privacy, reduces the communication effort, and can be implemented on a simple microcontroller. 
Future research should investigate the possible extension to non-convex storage, e.g. EVs with (dis-)charging efficiencies and restrictions on simultaneous (dis-)charging.

\bibliographystyle{IEEEtran}
\bibliography{root}

\end{document}